\documentclass[11pt,reqno]{amsart}
\usepackage{amsfonts}
\usepackage{mathrsfs}
\usepackage{amsmath}
\usepackage{amssymb,amsfonts}

\numberwithin{equation}{section}

\newcommand{\mysection}[1]{\section{#1}\setcounter{equation}{0}}

\newfont{\bb}{msbm10 at 12pt}

\def\R{\hbox{\bb R}}
\def\hH{\hbox{\bb H}}

\def\g{{\bar g}}

\def\k{{\lambda}}

\def\hconn{\breve{\nabla}}

\def\P{\mathcal P}

\def\T{\mathcal P}
\def\H{\mathcal H}

\def\hconn{\breve{\nabla}}

\newcommand{\bal}{\begin{aligned}}      \newcommand{\eal}{\end{aligned}}
\newcommand{\ba}{\begin{array}}      \newcommand{\ea}{\end{array}}
\newcommand{\bc}{\begin{center}}     \newcommand{\ec}{\end{center}}
\newcommand{\be}{\begin{enumerate}}  \newcommand{\ee}{\end{enumerate}}
\newcommand{\beq}{\begin{eqnarray}}  \newcommand{\eeq}{\end{eqnarray}}
\newcommand{\beQ}{\begin{eqnarray*}} \newcommand{\eeQ}{\end{eqnarray*}}
\newcommand{\bi}{\begin{itemize}}    \newcommand{\ei}{\end{itemize}}
\newcommand{\bt}{\begin{tabular}}    \newcommand{\et}{\end{tabular}}
\newcommand{\bdm}{\begin{displaymath}} \newcommand{\edm}{\end{displaymath}}

\newcommand{\ls}{\setlength{\baselineskip}{12pt}
                 \setlength{\parskip}{3mm}}

\begin{document}

\title[]{Recent progress on the positive energy theorem}

\author{Xiao Zhang}
\address[]{Institute of Mathematics, Academy of Mathematics and
Systems Science, Chinese Academy of Sciences, Beijing 100190, PR China}
\email{xzhang@amss.ac.cn}

\date{}

\begin{abstract}

We give a short review of recent progress on the positive energy theorem in general relativity, especially for
spacetimes with nonzero cosmological constant.

\end{abstract}

\maketitle \pagenumbering{arabic}

\mysection{Introduction}\ls

In general relativity, a spacetime is a 4-dimensional Lorentzian manifold $({\bf L}^{3,1}, \bf{g})$ which satisfies
the Einstein field equations
 \beQ
{\bf R} _{\alpha \beta}-\frac{{\bf R}}{2}{\bf g} _{\alpha \beta} +\Lambda {\bf g} _{\alpha \beta}
=T _{\alpha \beta},
 \eeQ
where ${\bf R} _{\alpha \beta}, {\bf R}$ are Ricci and scalar curvatures of ${\bf g}$ respectively,
$\Lambda $ is the cosmological constant and $T _{\alpha \beta}$ is the energy-momentum tensor of matter. Indeed,
$\Lambda $ is the Ricci curvature of vacuum spacetime $T _{\alpha \beta} =0$ where the Einstein field equations reduce to
$
{\bf R}_{\alpha \beta} =\Lambda {\bf g} _{\alpha \beta}.
$

The maximally symmetric vacuum solutions are as follows
\bi
\item [(i)] $\Lambda =0$, Minkowski spacetime $\R ^{3,1}$. The symmetric group is (3+1)-Poinca\'{r}e group generated by
four translational Killing vectors $L _\alpha$ and six rotational Killing vectors $U _{\alpha \beta}$.
\item[(ii)] $\Lambda >0$, de Sitter spacetime which is the hypersurface in $\R ^{4,1}$
\beQ
-\left(X^0 \right) ^2 +\left(X^1 \right) ^2 +\left(X^2 \right) ^2 +\left(X^3 \right) ^2 +\left(X^4 \right) ^2=\frac{3}{\Lambda}.
\eeQ
The symmetric group is (4+1)-Lorentzian group generated by ten rotational Killing vectors $U _{\alpha \beta}$.
\item [(iii)] $\Lambda <0$, anti-de Sitter spacetime which is the hypersurface in $\R ^{3,2}$
\beQ
-\left(X^0 \right) ^2 +\left(X^1 \right) ^2 +\left(X^2 \right) ^2 +\left(X^3 \right) ^2 -\left(X^4 \right) ^2=\frac{3}{\Lambda}.
\eeQ
The symmetric group is (3+2)-Lorentzian group generated by ten rotational Killing vectors $U _{\alpha \beta}$.
\ei

Choose a frame in spacetime ${\bf L}^{3,1}$ with $e_0$ timelike and $e_i$ spacelike, $T_{00}$ is interpreted as the local mass density and $T_{0i}$
is interpreted as the local momentum density. The dominant energy condition asserts that
 \beQ
T _{00} \geq \sqrt{T _{01} ^2 + T _{02} ^2 +T _{03} ^2}, \quad T _{00} \geq T _{\alpha \beta}.
 \eeQ
Let $(M, g, h)$ be an initial data set, where $M$ is 3-dimensional spacelike hypersurface in ${\bf L}^{3,1}$, $g$ is the Riemannian metric of
$M$ and $h$ is symmetric 2-tensor which is the second fundamental form of $M$. Let $\nabla$, $R$ be the Levi-Civita connection, the scalar curvature
of $M$ respectively. Then $g$ and $h$ are constrained by the Gauss and Codazzi equations which give
 \beQ
T_{00} =\frac{1}{2}(R +(h ^i _{\;\;i} ) ^2 -h _{ij} h ^{ij}),\quad T_{0i} = \nabla^j h _{ij} - \nabla _i h ^j _{\;\;j}.
 \eeQ

For spacetimes which are asymptotic to the these maximal symmetric solutions, the famous Noether theorem in the framework of a Lagrangian theory asserts
that the Killing fields of maximal symmetric spacetimes provide conserved quantities \cite{ADM, Chr1}. Physically, the dominant energy condition for the
energy-momentum tensor should indicate these conserved quantities are future timelike or null in suitable sense, and the above maximal symmetric
spacetimes have the lowest energy and serve as the ground states. But this physical predication is far from obvious due to the well-known pseudo-tensor
characteristic of the Lagrangian for gravity, e.g. \cite{Ch}. It was therefore conjectured to be true and refereed as the {\em positive energy conjecture}.
The conjecture plays a fundamental role and provides a self-consistent check for general relativity. However, what energy of the ground states is
needs much deeper physical theory.

The paper contributes to Proceedings of the International Conference on Gravitation and Cosmology/The Fourth Galileo-Xu Guangqi Meeting held on May 4-8
in Beijing. We would like to point out that many relevant interesting topics are not included due to limit pages of the paper.

\mysection{$\Lambda$ =0}\ls

In this section, we assume the cosmological constant is zero. We shall review the main total energy-momentum inequalities in this case.

\subsection{Positive energy theorem}

When the cosmological constant is zero and spacetimes are asymptotically flat, the time translation $L_0 $, the space translations $L_i $,
the space rotations $U _{ij} $ and the time-space rotations $U _{i0} $ give rise to the total energy, the total linear momentum,
the total angular momentum and the center of mass respectively.

An initial data set $(M, g, h)$ is {\it asymptotically flat} if there is a compact set $K \subset M$ such that $M \setminus K $ is the union
of finite $M _c$, where each $M_c \cong \R ^3 \setminus B _ r$ is called an ``end" of $M$, $B _r$ is the closed ball of radius $r$ with center
at the coordinate origin. In each end, let $\{x ^i \}$ be the Euclidean coordinates of $\R ^3$, $g_{ij}=\delta _{ij}+a_{ij}$, and $a_{ij}$ satisfy
 \beQ
a_{ij}=O\big(\frac{1}{r}\big), \partial _k  a_{ij}=O\big(\frac{1}{r ^2}\big), \partial _l \partial _k a _{ij}=O\big(\frac{1}{r ^3}\big),
h _{ij}=O\big(\frac{1}{r ^2}\big), \partial _k h _{ij}=O\big(\frac{1}{r ^3}\big),
 \eeQ
The total energy $E$, the total linear momentum $P_i$ for asymptotically flat initial data set $(M, g, h)$ were defined by Arnowitt, Deser and Misner \cite{ADM, ADM1}
 \beQ
E = \frac{1}{16\pi}\int _{S_{\infty }}\left(\partial _j g_{ij}-\partial _i g_{jj}\right)\ast dx ^i,\quad
P_{k} =\frac{1}{8\pi}\int _{S_{\infty}}\left(h_{ki}-g _{ki} tr_g(h)\right)\ast dx ^i
 \eeQ
where $S _{\infty}$ is the sphere at infinity on end $M_c$, $1 \leq k \leq 3$. Furthermore, we assume the absolute values of scalar curvature
$|R|$ and $|T_{0i}|$ are integrable over $M$. These yield that the total energy-momentum are geometric quantities which are independent on the
choice of coordinates $\{x ^i \}$ \cite{B, Chr1}.

The positive energy theorem for the ADM total energy-momentum was first proved by Schoen and Yau \cite{SY1,SY2,SY3}, then by Witten \cite{Wi1,PT}.
Schoen-Yau's proof included black hole's case and Witten's proof was extended to black hole's case in \cite{GHHP}.
Some higher dimensional cases were provided in \cite{EHLS, B, Di, XD} for instance.

{\em The positive energy theorem: If $(M, g, h)$ is asymptotically flat, with possibly a finite number of apparent horizons which are inner boundaries
of $M$, each of them is topological $S^2$ whose mean curvature $H$ satisfies $H \pm tr_g (h| _{S^2}) =0$. Suppose the dominant energy condition holds,
then\\
(i) $E \geq \sqrt{P _1 ^2 +P _2 ^2 + P_3 ^2}$ for each end;\\
(ii) That $E=0$ for some end implies $M$ has only one end, and ${\bf L} ^{3,1}$ is flat along $M$. If $h=0$, then $(M,g)$ is $\R ^3$.}

The total mass is $\sqrt{E^2 - P _1 ^2 -P _2 ^2 - P_3 ^2}$ which is a Lorentzian invariant. When the positive energy theorem holds, the total mass
is well-defined.

The positive energy theorem does not hold for spacetimes with naked singularity, e.g., Schwarzschild spacetime with $m<0$. So $E>0$ may be refereed as
the Schwarzschild constraint. The dominant energy condition provides sufficient conditions for the Schwarzschild constraint.

\subsection{Kerr Constraint}

Kerr spacetime also has naked singularity for $m<|a|$. So for general asymptotically flat spacetimes with the total energy $E$, the total angular
momentum $J$ and the total mass $m$, $E>|J|$ or $E>|J/m|$ for $m \neq 0$ may forbid naked singularities. Such an inequality is called
the {\em Kerr constraint} by Schoen recently.

In 1999, Zhang proved the positive energy theorem for {\em generalized asymptotically flat initial data set} $(M, g, p)$ where $p$ is general 2-tensor
which is not necessary symmetric, and satisfies the following conditions
 \beQ
a \in  C ^{2, \alpha} _{-\tau},\quad
 p \in  C ^{0, \alpha} _{-\tau -1},\quad tr_{\bar g}(p)
\in W ^{1,\frac{q}{2}} _{-\tau -1}, \quad \{d \theta , d ^* \theta \}\in  L _{\frac{q}{2}, -\tau -2}
 \eeQ
where $\theta =(p_{ij}-p_{ji})e^i \wedge e^j$ is the associated 2-form of $p$, $\frac{1}{2}< \tau \leq 1$ \cite{Z1}.

Geometrically, the second fundamental form $p$ is nonsymmetric when spacetimes equip with affine connections with torsion. In this case matter translates,
meanwhile, it rotates. The idea using connection with torsion was initially duo to E. Cartan \cite{Can, Can1, Can2, Can3}. In \cite{Z1}, Zhang defined the following
generalized linear momentum counting both translation and rotation
 \beQ
\bar P _{k} = \frac{1}{8\pi}\int_{S_{\infty}}(p_{ki}-g _{ki} tr_g (p))\ast dx ^i,
\eeQ
and proved the following theorem.

{\em The generalized positive energy theorem: If $(M, g, p)$ is generalized asymptotically flat, with a finite number of generalized apparent horizons
which are inner boundaries of $M$, each of them is topological $S^2$ whose mean curvature $H$ satisfies $H \pm tr_g (p| _{S^2}) =0$. Suppose the
generalized dominant energy condition
 \beQ
\frac{1}{2}\big(R +(tr _g (p))^2 -|p|^2\big)\geq\max \big\{|\omega|, |\omega+\chi|\big\},
 \eeQ
holds, $\omega _j = \nabla ^i p _{ji}-\nabla _j tr_g(p)$, $\chi _j =2 \nabla ^i (p _{ij} - p _{ji})$ (the term in left hand side can be interpreted as
$T_{00}$ and the two terms in right hand side can be interpreted as $|T_{0i}|$, $|T_{i0}|$ in suitable sense), then\\
(i) $E \geq \sqrt{\bar P_1 ^2 + \bar P _2 ^2 + \bar P _3 ^2}$ for each end;\\
(ii) That $E=0$ for some end implies $M$ has only one end, and
\beQ
 R_{ijkl}+p _{ik}p _{jl}-p _{il}p _{jk} = 0,\quad
 \nabla _i p _{jk} -\nabla _j p _{ik} = 0,\quad
 \nabla ^i (p _{ij} -p _{ji}) = 0.
\eeQ}

In Newton mechanics, let $\{x_0 ^1, x_0 ^2, x_0 ^3\}$ be the system's center of mass, $T_{v}$ be the momentum density of the system, the total angular
momentum of rigid body ${\bf V}$ is
 \beQ
 J_{k} = \int _{\bf V}\epsilon _{kuv} (x ^u -x_0 ^u ) T ^{v} *1
       = \int _{S _\infty}\epsilon _{kuv} (x ^u -x_0 ^u ) (h ^v _{\;\;i} -\delta _{vi} tr (h)) *dx^i.
 \eeQ
The second equality holds if, furthermore, ${\bf V}=\R ^3$, $T _{v}=\partial _i h _v ^{\;\;i} -\partial _v tr(h)$. Usually, the total angular momentum
is a set of anti-symmetric two tensors $J_{ij}=-J_{ji}$. In $\R ^3$, the rotation along a plane is equivalent to the
translation alone an axis which is perpendicular to the plane, $J _k = \epsilon _{k} ^{\;\;ij}J _{ij}$. So we can define them as co-vectors in
3-dimensional space. But the above definition of $J_k$ does not make sense in higher dimension.

In 1974, Regge-Teitelboim generalized it to general relativity and defined the total angular momentum for asymptotically flat initial data sets \cite{RT},
 \beQ
 J _k (x_0) =\frac{1}{8\pi}\int_{S_{\infty}} \epsilon _{kuv} (x ^u -x_0 ^u )\pi ^v _{\;i} *dx ^i, \quad
 \pi ^v _{\;i} =h ^{v} _{\;i}-g ^{v} _{\;i} tr_g (h).
 \eeQ
In general, the integrand is O($\frac{1}{r}$) which may not be integrable. This ambiguity resolution requires stronger ``Regge-Teitelboim"
conditions on ends
\beQ
g(x)-g(-x)=O(r^{-3}), \quad \pi (x) +\pi (-x) =O(r^{-3}).
\eeQ

However, as the integrand is not tensor, it is not possible to relate the local density to Regge-Teitelboim's total angular momentum. To resolve
this difficulty, Zhang defined trace free, non-symmetric tensor of local angular momentum density \cite{Z1}
\beQ
\tilde h ^z _{ij} = \frac{1}{2}\epsilon _{i} ^{\;\;\;uv}\big(\nabla _u \rho _z ^2 \big)\big(h _{vj}- g _{vj} tr _g (h)\big)
\eeQ
where $\rho _z$ is the distance function w.r.t some $z \in M$. If $(M,g, \tilde h ^z _{ij})$ is generalized asymptotically flat,
Zhang defined the total angular momentum \cite{Z1}
\beQ
J_{k} =\frac{1}{8\pi}\int _{S_{\infty}}\tilde h ^z _{ki}\ast dx ^i.
\eeQ
$J$ is also a geometric quantity which is independent on the choice of coordinates $\{x^i\}$ if $|\nabla \tilde h ^z|$ is integrable.
In Kerr spacetime, $h_{ij}=O(\frac{1}{r^4})$, so the total angular momentum is well-defined and it was found $J=(0, 0, ma)$ \cite{Z1-1}.
By taking $p_{ij}=C \tilde h ^z _{ij}$ for certain constant $C>0$, Zhang proved the Kerr constraint under the generalized dominant energy
condition \cite{Z1}. The higher dimensional Kerr constraint is interesting in mathematics and it is still open.

The dominant energy condition does not yield to the Kerr constraint. In \cite{HSW}, Huang, Schoen and Wang showed that it is possible
to perturb arbitrary vacuum asymptotically flat initial data sets to new vacuum ones having exactly the same total energy, but with the
arbitrary large total angular momentum. However, the Kerr constraint holds for simply connected, asymptotically flat, maximal, axisymmetric,
vacuum black hole initial data sets $(M, g, h)$ defined on $\R ^3 \setminus \{\rho =0\}$
 \beQ
g = e^{2\alpha -2U}\big(d\rho ^2 +dz^2 \big)+e^{-2U} \rho^2 \big(d\phi +\rho B d\rho +A dz\big)^2, \quad tr_g(h)= 0,
 \eeQ
where $(\rho, \phi, z)$ are cylindrical coordinates, and all functions are $\phi$ independent \cite{D, Chr2, CLW, SZ}.

\subsection{Bondi energy-momentum}

Bondi-Sachs' radiating metrics for gravitational waves are wave-like, vacuum solutions of the Einstein field equations,
and take the following asymptotical forms
 \beQ
 \begin{aligned}
{\bf g} _{BS} =&-\big(1-\frac{2M}{r}\big)du ^2 -2du dr+2l du
d\theta +2 \bar l \sin \theta du d\psi \\
 &+r ^2 \Big[(1+\frac{2c}{r})d\theta ^2 +\frac{4d}{r}\sin \theta d
\theta d \psi
+(1-\frac{2c}{r})\sin ^2 \theta d \psi ^2\Big]\\
 &+\mbox{lower order terms}
 \end{aligned}
 \eeQ
where $u$ is retarded coordinate ($u$-slices are null hypersurfaces), $r=x^1$, $\theta $ and $\psi$ are polar coordinates,
$M$, $c$, $d$ are smooth functions of $u$, $\theta$, $\psi$ defined on $\R \times S ^2$ with regularity condition
$\int _0 ^{2\pi} c(u, \theta, \psi)d\psi =0$ for $\theta =0, \pi$ and all $u$, $l = c _{, \theta} +2c \cot \theta +d _{, \psi} \csc \theta$,
$\bar l = d_{, \theta} +2d \cot \theta -c _{,\psi} \csc \theta$.

Throughout the paper, we denote $n^0=1$, $n^1 =\sin \theta \cos \psi$, $n^2 =\sin \theta \sin \psi$, $n^3=\cos \theta$. At null infinity,
the Bondi energy-momentum of $u _0$-slice are defined as
\begin{eqnarray*}
m _\nu (u _0) = \frac{1}{4 \pi} \int _{S ^2} M (u _0, \theta, \psi)n ^{\nu} d S
 \end{eqnarray*}
for $\nu =0,1,2,3$. The famous Bondi energy loss formula asserts
 \beQ
\frac{d}{du} m _{0} (u)=-\frac{1}{4\pi}\int _{S ^2} (c _{,u}) ^2 +(d _{,u})^2 \leq 0.
 \eeQ
Let $|m(u)|=\sqrt{m _1 ^2 (u)+ m _2 ^2 (u)+m _3 ^2 (u)}$. If $|m(u)|\neq 0$, it was generalized to the Bondi energy-momentum loss \cite{HYZ}
 \beQ
\frac{d}{du} \Big(m _0 (u) -|m(u)| \Big) = -\frac{1}{4\pi}\int _{S ^2} \Big[(c _{,u}) ^2 +(d _{,u})^2 \Big]
\Big(1- \frac{m_i n^i}{|m|} \Big) \leq 0.
 \eeQ
This asserts that the Bondi energy-momentum can be viewed as the total energy-momentum measured after the loss due to the gravitational
radiation up to that time.

It is an old question whether isolated gravitational systems can radiate away more energy than they have initially, i.e., whether Bondi's Energy
should be nonnegative. The proofs of this positivity were claimed by using both Schoen-Yau and Witten's positive energy arguments (see \cite{CJL} and
references therein). However, extra conditions are required when two methods are worked out rigorously and completely \cite{HYZ}.

Let ${\mathcal M}(u, \theta, \psi)=M(u, \theta, \psi) -\frac12\big(l _{,\theta} +l \cot \theta+\bar l _{,\psi} \csc \theta \big)$. In \cite{HYZ}, it was
derived that ${\mathcal M} _{,u} =-c _{,u} ^2 -d _{,u} ^2$. Suppose $m_0 (u)=|m(u)|$ on $[u_1, u_0]$. If $m_0 (u) \neq 0$, then the Bondi
energy-momentum loss formula implies $m _i (u)=|m(u)| n ^i$. This is impossible except $m _{i}=0$ since $m _{i}$ are independent on $\theta$, $\psi$
but $n ^i$ do depend on $\theta$, $\psi$. Thus $m_0 (u) = 0$ on $[u_1, u_0]$. By the Bondi energy loss formula, we obtain that
$c _{,u}= d _{,u}={\mathcal M} _{,u} =0$ on $[u_1, u_0]$. Therefore,\\\\
(i) if ${\mathcal M} (u_0, \theta, \psi)=C$ for some constant $C$, then ${\mathcal M}(u, \theta, \psi)=C$ on $[u_1, u_0]$;\\
(ii) if $c(u_0, \theta, \psi) =d (u_0, \theta, \psi)=0$, then $c(u, \theta, \psi) =d(u, \theta, \psi) =0$ on $[u_1, u_0]$.

Now we discuss the positivity of Bondi energy-momentum. The following two theorems are proved in \cite{HYZ} and stated here more accurately.
Firstly, Schoen-Yau's method shows

{\em Positivity of Bondi energy (Schoen-Yau's method): Suppose there exists $u _0$ in vacuum Bondi's radiating spacetime such that
${\mathcal M}(u _0, \theta, \psi)$ is constant.\\
(i) $m _{0} (u_0) \geq |m(u_0)|$, and the Bondi energy-momentum loss formula gives $m _{0} (u) \geq |m(u)|$
for all $u \leq u_0$;\\
(ii) If $m _{0} (u_0)=|m(u_0)|$ and there is $u_1 < u_0$ such that $m _{0} (u_1)=|m(u_1)|$, then ${\mathcal M}$ is constant on $[u_1, u_0]$.
Thus the spacetime is flat in the region $(u_1, u_0]$.}

The proof of the positivity using Witten's method requires the positive energy theorem for asymptotically null initial data sets.
In Minkowski spacetime $\R^{3,1}$, the null cone is $\{t=r\}$. The spacelike hypersurface $t=\sqrt{1+r ^2}$ approaches to null infinity, and equips with
the hyperbolic metric $\breve{g}$ and the nontrivial second fundamental form $\breve{h}$ in polar coordinates
 \begin{eqnarray*}
 \breve{g} = \breve{h} =\frac{ dr ^2 }{1+r ^2}+r ^2 \big(d \theta ^2 +
\sin ^2 \theta d\psi ^2\big).
 \end{eqnarray*}
Denote $\breve{e} _i$, $\breve{e}^i$ and $\hconn $ the frame, coframe and the Levi-Civita connection respectively. The initial data set $(M, g, h)$ is
asymptotically null if, on each end,
 \beQ
\Big\{a_{ij},\,\, \hconn  _k a_{ij},\,\, \hconn  _l \hconn  _k a_{ij},\,\,
b _{ij}, \,\,\hconn  _k b_{ij} \Big\} = O(r^{-\tau})
 \eeQ
for $\tau >\frac{3}{2}$, where $a_{ij}=g(\breve{e} _i,\breve{e} _j)-\breve{g}(\breve{e}_i,\breve{e}_j)$,
$b_{ij}=h (\breve{e} _i,\breve{e} _j )-\breve{h}(\breve{e}_i,\breve{e}_j)$. Let $R$, $\nabla$, $\rho_z$ be the scalar curvature and the Levi-Civita
connection of $g$, and the distance function with respect to $z \in M$ respectively. We further assume
\beQ
(R +6)\rho_z \in L^1(M),\quad \big(\nabla^j  h _{ij}-\nabla_i tr_{g}(h)\big)\rho_z \in L^1(M).
\eeQ
The asymptotic null total energy-momentum of the end $M _l$ are
 \beQ
E _{\nu}=\frac{1}{16\pi}\int_{S_{\infty}}\mathcal{E} n^\nu r \breve{e}^2\wedge \breve{e}^3,
 \eeQ
where $\nu =0,1,2,3$, $\mathcal{E} =\hconn ^{j}  g_{1j}- \hconn  _1 tr_{\breve{ g} }(g ) +(a_{22}+a_{33})+2(b_{22}+b _{33})$.

{\em The positive energy theorem near null infinity \cite{Z2}: Let $(M, g, h)$ be an asymptotically null initial data set in spacetime
${\bf{L}}^{3,1}$. Suppose ${\bf{L}}^{3,1}$ satisfies the dominant energy condition, then\\
(i) $E_{0} \geq \sqrt{ E_{1}^2 +E_{2}^2 +E_{3}^2}$ for any end;\\
(ii) That $E_{0}=0$ for some end implies $M$ has only one end, and ${\bf L} ^{3,1}$ is flat along $M$.}

It ensures that the asymptotically null total mass $\sqrt{E_0 ^2 -E_1 ^2 -E_2 ^2 -E_3 ^2}$ is well-defined.

In order to using this theorem to provide the positivity of Bondi energy, we use the spacelike hypersurface
\beQ
u=u_0 +\sqrt{1+r^2} -r +\frac{c^2(u_0, \theta, \psi) +d^2(u_0, \theta, \psi)}{12r^3}+\frac{a_3(\theta, \psi)}{r^4}+o(\frac{1}{r^4})
\eeQ
to cut the Bondi-Sachs metrics and obtain asymptotically null initial data sets. Unfortunately, the metric and the second fundamental form of this
spacelike hypersurface differ from the hyperbolic metric with the leading terms $\frac{c(u_0)}{r}$ and $\frac{c(u_0)}{r}$ \cite{HYZ}. This breaks the
conditions required in the above positive energy theorem near null infinity. So we need to assume $c(u_0, \theta, \psi)=d(u_0, \theta, \psi)=0$. Under
these conditions
\beQ
E_\nu =\frac{1}{16\pi}\int _{S^2 _\infty} \mathcal{E} n^\nu r \breve{e}^2 \wedge \breve{e}^3= m_\nu (u_0)
\eeQ
on the asymptotically null initial data sets. (The notations here are different from those in \cite{HYZ}, $E_\nu$, $\mathcal{E}$ here are
$E_\nu (\mathbb{X}) -P_\nu (\mathbb{X})$, $\mathcal{E} -2\mathcal{P}$ in \cite{HYZ} respectively. It is an error to integrate
$\mathcal{E} -\mathcal{P}$ to get $E_\nu (\mathbb{X}) -P_\nu (\mathbb{X})=\frac{5}{8}m_\nu (u_0)$ in \cite{HYZ} which should be corrected as above.)
Then the positive energy theorem near null infinity shows \cite{HYZ}

{\em Positivity of Bondi energy (Witten's method): Suppose there exists $u _0$ in vacuum Bondi's radiating spacetime such that $c(u_0)=d(u_0)=0$.\\
(i) $m _{0} (u_0) \geq |m(u_0)|$, and the Bondi energy-momentum loss formula gives $m _{0} (u) \geq |m(u)|$
for all $u \leq u_0$;\\
(ii) If $m _{0} (u_0)=|m(u_0)|$ and there is $u_1 < u_0$ such that $m _{0} (u_1)=|m(u_1)|$, then $c(u_0)=d(u_0)=0$ on $[u_1, u_0]$. Thus the spacetime
is flat in the region $(u_1, u_0]$.}

It ensures that the Bondi mass $\sqrt{m_0 ^2 -m_{1} ^2 -m _{2} ^2 -m_{3} ^2}$ is well-defined.

It is an interesting question whether the total angular momentum can be detected near or at null infinity. As it is still open to find smooth
Bondi-Sachs' coordinates for Kerr spacetime, it is unclear what rotation means for gravitational systems traveling in the speed of light. On the other hand,
B. Liu constructed certain asymptotically null initial data sets in Kerr spacetime on the dissertation of his Bachelor's degree at the University of
Science and Technology of China in 2013. Choose the spacelike hypersurface
$$
t=r+2m\ln r +\frac{c_1(\theta, \psi)}{r}+\frac{c_2(\theta, \psi)}{r^2}+\frac{c_3(\theta, \psi)}{r^3}+\cdots
$$
in Kerr spacetime. Liu found the fall-offs of the metric and the second fundamental form by suitable choices of $c_i(\theta, \psi)$
\beQ
\begin{aligned}
g_{11}=&1+O(\frac{1}{r^4}), \quad g_{12}=\frac{2 m^2 a^2 \sin 2\theta}{3r^3}+O(\frac{1}{r^4}),\\
g_{13}=&-\frac{2m a \sin \theta}{r} -\frac{4 m^2 a \sin \theta}{r^2}-\frac{ma(8m^2-5a^2\cos^2 \theta +a^2)\sin \theta}{r^3}+O(\frac{1}{r^4}),\\
g_{22}=&1+O(\frac{1}{r^5}), \quad g_{23}=O(\frac{1}{r^5}), \quad g_{33}=1+\frac{2m a^2 \sin^2\theta}{r^3}+O(\frac{1}{r^5}),\\
h_{11}=&1-\frac{m^2 a^2 \sin ^2 \theta}{r^2} -\frac{m(16m^2 a^2-3a^2)\sin^2 \theta +2m}{r^3}+O(\frac{1}{r^4}),\\
h_{12}=&\frac{8m^2 a^2 \sin 2 \theta}{3r^3}+O(\frac{1}{r^4}),\quad h_{23}=O(\frac{1}{r^4}),\\
h_{13}=&-\frac{2 m a \sin \theta}{r} -\frac{4 m^2 a \sin \theta}{r^2}\\&
        -\frac{ma\big[4a^2(4m^2+5)\sin^2 \theta +16m^2 -4a^2 +3\big]\sin \theta}{r^3}+O(\frac{1}{r^4}),\\
h_{22}=&1+\frac{2m^2 a^2 \sin ^2 \theta}{r^2} -\frac{m(4m^2 a^2 \sin^2 \theta -1)}{r^3}+O(\frac{1}{r^4}),\\
h_{33}=&1+\frac{2m^2 a^2 \sin ^2 \theta}{r^2} +\frac{m\big[(8m^2-1)a^2 \sin^2 \theta +1\big]}{r^3}+O(\frac{1}{r^4}).
\end{aligned}
\eeQ

Let $\tilde h _{ij}=h_{ij}-g_{ij}$, $\tilde {\pi} _{ij} =\tilde h_{ij} -tr _g(\tilde h)$. If $a \neq 0$, Liu found,
\beQ
\tilde \pi _{12}=O(\frac{1}{r^3}), \quad \tilde \pi _{13}=O(\frac{1}{r^3}), \quad \mbox{other} \,\,\tilde \pi _{ij}=O(\frac{1}{r^4}).
\eeQ
The fall-offs of $g_{13}$, $\tilde \pi _{12}$ and $\tilde \pi _{13}$ are not sufficiently fast to define the total angular momentum.
This indicates the rotation does not make sense also near null infinity.

\mysection{$\Lambda$ >0}\ls

In this section, we assume the cosmological constant is positive. As recent cosmological observations indicated that our universe should
have a positive cosmological constant, the positive energy theorem for asymptotically de Sitter spacetimes gains much more importance. There are a large
number of papers contributing to issue of the total energy-momentum in spacetimes with the positive cosmological constant, e.g.,
\cite{AD, BBM, CMZ, AMR, GM, ABK1, ABK2}, but only few studying their positivity \cite{S, SIT, KT}. Even through, the formulations of the
positivity are incomplete, and some of the proofs do not seem correct in mathematics. In \cite{LXZ, LZ}, the complete and rigorous study of the
positive energy theorem were provided.

Denote $\lambda =\sqrt{\frac{3}{\Lambda}}$ throughout the section. De Sitter spacetime can be fully covered by the global coordinates equipped with the
de Sitter metric
\beQ
\tilde g _{dS} = -d \bar{t} ^2 + \k^2 \cosh^2 \frac{\bar{t}}{\k}
\Big(d \bar{r} ^2 +\sin ^2 \bar{r} (d \bar \theta ^2 +\sin ^2 \bar \theta d \bar \psi ^2) \Big).
\eeQ
In global coordinates, each time slice is a 3-sphere with constant curvature and has no spatial infinity. Therefore the ADM formulation of the
energy-momentum is not available. There are essentially two ways to separate de Sitter spacetime into two parts which give two different spatial
infinities. And the positive energy theorems can be established for spacetimes which are asymptotic to either half of the de Sitter spacetime under
reasonable conditions.

Separating along the hypersurface $X^0=X^4$, the half-de Sitter spacetime is covered by planar coordinates equipped with the de Sitter metric
\beQ
\tilde g _{dS} = -dt ^2 + e ^\frac{2t}{\k} g_\delta , \quad g_\delta=(dx ^1)^2 +(dx ^2)^2+(dx ^3)^2
\eeQ
and the initial data set is $(\R ^3, \breve{g}= e ^\frac{2t}{\k} g_{\delta}, \breve{K}=\frac{1}{\k} \breve{g})$.
Separating along the hypersurface $X^4=-\k$ ($X^4=\k$), the half-de Sitter spacetime $X^4 <-\k$ ($X^4 >\k$) is covered by hyperbolic coordinates
equipped with the de Sitter metric
 \beQ
\tilde g _{dS} = -dT^2 + \sinh^2 \frac{T}{\k}\breve{g} _{\H}, \quad \breve{g} _{\H}=dr^2+\k^2\sinh^2 \frac{r}{\k}(d\theta^2+\sin^2\theta d \psi^2)
 \eeQ
and the initial data set is $(\hH ^3, \breve{g}=\sinh^2 \frac{T}{\k}\breve{g} _{\H}, \breve{K}=\frac{1}{\k}\coth \frac{T}{\k} \breve{g})$ where $\breve{g} _{\H}$
is the standard hyperbolic metric.

To define the total energy-momentum from the Hamiltonian point of view, an asymptotically de Sitter initial data set $(M, g, K)$ should satisfy
\beQ
g-\breve{g}=e ^\frac{2 t_0}{\k} a =O(\frac{1}{r}), \qquad K-\breve{K}=e ^\frac{t_0}{\k} b =O(\frac{1}{r^2})
\eeQ
for constant $t_0$ on ends in planar coordinates. Then the ten Killing vectors $U_{\alpha \beta}$ of $\R^{4,1}$ and $a$ and $b$ are used to define the
ADM-like total energy-momentum. But there is no energy-momentum inequalities for them in general. Alternatively, there is different approach to define the
total energy-momentum from the initial data set point of view. We explain the main idea for this approach in planar coordinates for simplicity, and it is similar
in hyperbolic coordinates.

In general, an initial data set $(M, g, K)$ in asymptotically de Sitter spacetimes should take the forms $g=e ^{2u(x)} \bar g (x)$, $h=K-\frac{1}{\k} g=e ^{u(x)} \bar h (x)$,
$u(x) = \frac{t_0}{\k}+o(1)$ and $(M, \bar g, \bar h)$ is asymptotically flat. In order to preserve the asymptotic flatness for $g$ and $h$, one must take
$u (x) = \frac{t_0}{\k} +\frac{a_1}{r}+O(\frac{1}{r^2})$. In this case it can be written that $g= e^{\frac{2 t_0}{\k}} \bar g$,  $h=e ^{\frac{t_0}{\k} } \bar h $
on ends, and the ADM total energy-momentum of $\bar g$, $\bar h$ can serve as the total energy-momentum of $g$, $h$ up to certain constant factors.
As $e^{\frac{t_0}{\k}}$ is constant, it is equivalent to assume $g= e^{\frac{2 t_0}{\k}} \bar g$, $h=e ^{\frac{t_0}{\k} } \bar h $ on whole $M$, instead of on ends only.

Thus, an initial data set $(M, g, K)$ is {\em $\P$-asymptotically de Sitter} if $g= \P ^2 \bar g$, $h= \P \bar h$ for certain constant
$\P>0$, and $(M, \bar g, \bar h)$ is asymptotically flat. Let $\bar E$, $\bar P _{k}$ and $\bar {J _{k}} (z)$ be the total energy, the total linear
momentum and the total angular momentum of the end $M_l$ for $(M,\bar g, \bar h)$ respectively. The corresponding quantities for $\P$-asymptotically
de Sitter initial data set $(M, g, K)$ are
\beQ
E=\P \bar E,\quad P_{k} =\P^2 \bar P _{k}, \quad J _{k}(z)=\P^2 \bar {J _{k}} (z).
\eeQ
and the positive energy theorem proved by Luo, Xie and Zhang \cite{LXZ} is as follows.

{\em The positive energy theorem (planar coordinates): Let $(M, g, K)$ be a $\P$-asymptotically de Sitter initial data set in spacetime
${\bf{L}} ^{3,1}$ with positive cosmological constant $\Lambda >0$. Suppose
\beQ
tr _g (K) \leq \sqrt{3\Lambda }
\eeQ
If ${\bf{L}}^{3,1}$ satisfies the dominant energy condition, then\\
(i) $E \geq \sqrt{P_{1} ^2 +P_{2} ^2 +P_{3} ^2}$ for any end;\\
(ii) That $E=0$ for some end implies
\beQ
(M, g, K) \equiv \Big(\R ^3, \P ^2 g_{\delta}, \sqrt{\frac{\Lambda}{3}}\P ^2 g_{\delta} \Big)
\eeQ
and the spacetime ${\bf{L}}^{3,1}$ is de Sitter along $M$.}

It ensures that the total mass $\sqrt{E^2 -P_1 ^2 -P_2 ^2 -P_3 ^2}$ is well-defined.

The total angular momentum $J$ is computed to equal $(0,0,\frac{ma}{(1+\frac{a^2}{\lambda ^2})^2})$ for Kerr-de Sitter spacetime \cite{LXZ}. It is
interesting that $J _3(z)$ conjugates to $J_{23}$ of \cite{HT} by replacing the positive cosmological constant to the negative cosmological constant.
The corresponding Kerr constraint was also proved in \cite{LXZ}.

{\em Let $(M, g, K)$ be $\T $-asymptotically de Sitter which has no apparent horizon in spacetime $(N^{1,3},\tilde{g})$ with positive cosmological
constant $\Lambda >0$. Denote $p=C \tilde {{\bar h} } ^z $ where $\tilde {{\bar h}} ^z $ is local angular momentum density defined by $\bar h$,
and $C>0$ is certain constant. Suppose that there exists a point $z \in M$ such that $(M, \bar{g}, p )$ is generalized asymptotically flat. If
$(M, \bar{g}, p)$ satisfies the generalized dominant energy condition, then\\
(i) $E \geq C |J(z)|_{\breve{g} _{\T}}$ for any end;\\
(ii) That $E=0$ for some end implies that $M$ has only one end, and
\beQ
 R_{ijkl}+p _{ik}p _{jl}-p _{il}p _{jk} = 0,\quad
 \nabla _i p _{jk} -\nabla _j p _{ik} = 0,\quad
 \nabla ^i (p _{ij} -p _{ji}) = 0.
\eeQ}

In the Hamiltonian formulation, $g-\breve{g}=O(r^{-1})$ and $K-\breve{K}=O(r^{-2})$ are used to define the total energy-momentum. In the initial data set
formulation, $g-\breve{g}=O(r^{-1})$, $h=O(r^{-2})$ which give $K-\breve{K}=O(r^{-1})$. Although the two approaches give the same total energy, the total
momenta are completely different. While it is finite, defined in terms of $h$, the total linear momentum is infinite in general, defined in terms of
$K-\breve{K}$. So there should not be any energy-momentum inequality in the Hamiltonian formulation. But it is known that $E\geq 0$ if
$tr _g (K) \leq \sqrt{3\Lambda }$ in this case.

In hyperbolic coordinates, denote $h=K-\frac{1}{\k}\coth \frac{T}{\k} g$, $\H =\sinh \frac{T}{\k}$ for certain constant $T$. Let $\breve{e} _i$,
$\breve{e}^i$ and $\hconn ^\H$ be the frame, coframe and the Levi-Civita connection of the hyperbolic metric $\breve{g} _{\H}$. An initial data set
$(M, g, K)$ is {\em $\H$-asymptotically de Sitter} if $g= \H ^2 \bar g$, $h= \H \bar h$, and, on each end, $\bar g$ and $\bar h$ satisfy
 \beQ
\Big\{a_{ij},\, \hconn ^\H _k a_{ij},\, \hconn ^\H _l \hconn ^\H _k a_{ij},\,\bar h _{ij}, \,\hconn ^\H _k \bar h_{ij} \Big\}=O(e^{-\frac{\tau}{\k}r})
 \eeQ
for $\tau >\frac{3}{2}$, $a_{ij}=\bar g(\breve{e} _i,\breve{e} _j)-\breve{g}_\H(\breve{e}_i,\breve{e}_j)$,
$\bar h_{ij}=\bar h (\breve{e} _i,\breve{e} _j )$. Let $\bar R$, $\bar \nabla$, $\rho_z$ be the scalar curvature and the Levi-Civita connection of
$\bar g$, and the distance function with respect to $z \in M$ respectively. We further assume
 \beQ
(\bar R +\frac{6}{\k^2})e^{\frac{\rho_z}{\lambda}},\quad
(\bar \nabla^j \bar h_{ij}-\bar \nabla_i tr_{\bar g}(\bar h))e^{\frac{\rho_z}{\lambda}} \in L^1(M).
 \eeQ
The total energy-momentum of the end $M _l$ are
 \beQ
E^{\H} _{\nu}=\frac{\H ^2}{16\pi}\int_{S_{\infty}}\mathcal{E} n^\nu e^{\frac{r}{\k}} \breve{e}^2\wedge \breve{e}^3,
 \eeQ
where $\nu =0,1,2,3$, $\mathcal{E} =\hconn ^{\H,j} \bar g_{1j}- \hconn ^{\H} _1 tr_{\breve{ g} _\H}(\bar g )
+\frac{1}{\k}(a_{22}+a_{33})+2(\bar h_{22}+\bar h _{33})$. The positive energy theorem in this case was also proved by Luo, Xie and Zhang \cite{LXZ}.

{\em The positive energy theorem (hyperbolic coordinates): Let $(M, g, K)$ be a $\H$-asymptotically de Sitter initial data set
in spacetime ${\bf{L}} ^{3,1}$ with positive cosmological constant $\Lambda >0$. Suppose
 \beQ
tr _g (K) \sinh\frac{T}{\lambda} \leq \sqrt{3\Lambda }\cosh\frac{T}{\lambda}.
 \eeQ
If ${\bf{L}} ^{3,1}$ satisfies the dominant energy condition, then\\
(i) $E ^{\H} _{0} \geq \sqrt{ (E ^{\H}) _{1} ^2 +(E ^{\H}) _{2} ^2 +(E ^{\H}) _{3} ^2 }$ for each end;\\
(ii) That $E^{\H} _0 =0$ for some end implies
\beQ
(M, g, K) \equiv \Big(\hbox{\bb H} ^3, \sinh^2 \frac{T}{\lambda} \breve{g}_{\H}, \sqrt{\frac{\Lambda}{3}}
\sinh\frac{T}{\lambda}\cosh\frac{T}{\lambda}\breve{g}_{\H} \Big)
\eeQ
and the spacetime ${\bf{L}}^{3,1}$ is de Sitter along $M$.}

It ensures that the total mass $\sqrt{ (E^{\H} _0 )^2 - (E ^{\H}) _{1} ^2 -(E ^{\H}) _{2} ^2 -(E ^{\H}) _{3} ^2 }$
is well-defined.

The positive energy theorems hold also when $M$ has a finite number of inner boundaries, and each $(\Sigma,\bar{g},\bar{h})$ of them is
topological $S^2$ whose induced metric $\bar{g}$ and the second fundamental form $\bar{h}$ satisfy
\beQ
\pm tr_{\bar g} (h|_\Sigma ) =tr _{g _\Sigma }(K \big| _\Sigma )-2 \sqrt{\frac{\Lambda }{3}}
\eeQ
in $\T $-asymptotically de Sitter initial data sets and
\beQ
\pm tr_{\bar g} (h|_\Sigma ) =tr _{g _\Sigma }(K \big| _\Sigma )-2\sqrt{\frac{\Lambda }{3}} \tanh\frac{T}{2\lambda}
\eeQ
in $\H $-asymptotically de Sitter initial data sets. All these inner boundaries do not coincide with the future/past apparent horizon
\beQ
\pm tr_{\bar g} (h|_\Sigma ) =tr _{g _\Sigma }(K \big| _\Sigma )
\eeQ
defined by the outward expansion of the future and past going light rays emanating from $\Sigma$.

However, as pointed out by Witten \cite{Wi2}, there is no positive conserved energy in de Sitter spacetime, and the corresponding Killing vector fields
to the Lorentzian generators are timelike in some region of de Sitter spacetime and spacelike in some other region (see also \cite{AD}). This indicates
that there should not have the positive energy theorem in the standard sense without extra mean curvature constraints.

The positive energy theorem indicates that any metric on $\R^n$ ($n\geq 3$) with scalar curvature $\geq 0$ must be flat if it agrees with the Euclidean
metric outside a compact set. In 1989, Min-Oo proved a theorem which implies any metric on $\hH^n$ with scalar curvature $R \geq -n(n-1)$ ($n\geq 3$) must
be hyperbolic if it agrees with the hyperbolic metric outside a compact set \cite{Mi}. In global coordinates, the de Sitter spacetime can be separated
into two parts along the hypersurface $X^3=0$. The time slices in this half-de Sitter spacetime are hemispheres. In 1995, Min-Oo conjectured that any
metric on the hemisphere $S^n_+$ ($n\geq 3$) with scalar curvature $R \geq n(n-1)$ must be the standard metric on $S^n_+$ if the boundary $\partial S^n_+$
is totally geodesic and the induced metric on $\partial S^n_+$ agrees with the standard metric on $\partial S^n_+$. In 2010, Brendle, Marques and Neves
constructed some metrics on the hemisphere $S^n_+$ ($n\geq 3$) with scalar curvature $R\geq n(n-1)$ and $R>n(n-1)$ at some point, but agree with the
standard metric in a neighborhood of $\partial S^n_+$ \cite{BMN}. It therefore provides a counterexample to Min-Oo's conjecture and provide one failure
of the positive energy theorem in certain case for positive cosmological constant.

In 2012, Liang and Zhang constructed certain asymptotically de Sitter initial data sets with negative total energy in de Sitter spacetimes \cite{LZ}.
These initial data sets violate the mean curvature constraints in two positive energy theorems.

In planar coordinates, consider the graph $f(x)=t_0+\varepsilon (1+r^2 )^{-\frac{1}{2}}$ for certain constant $t_0$ in de Sitter spacetime. The graph is
$\P$-asymptotically de Sitter ($\P=e^{\frac{t_0}{\k}}$) for sufficiently small $\varepsilon$, and
 \beQ
 \begin{aligned}
g_{ij}=&e^{\frac{2f}{\k}}\delta _{ij} - \varepsilon ^2 \frac{x_i x_j}{(1+r^2)^3} ,\\
tr _g (K)=&\frac{3}{\k} + \frac{1}{2\k}\varepsilon ^2 \frac{r^2}{ (1+r^2)^3}e^{-\frac{2f}{\k}} + O(r^{-5}).
 \end{aligned}
 \eeQ
Thus if $ \varepsilon \neq 0$, then $tr _g (K) > \frac{3}{\k}$ for large $r$. Choose $\varepsilon <0$, the total energy
 \beQ
 \begin{aligned}
    E &=\frac{\P}{16\pi}\lim _{r \rightarrow \infty} \int
         _{S_r}(\partial _j \g _{ij}-\partial _i \g _{jj})*dx ^i\\
      &=\frac{\P}{16\pi}\lim _{r \rightarrow \infty} \int _{S_r}\Big(
           \frac{4\varepsilon }{\k}\frac{r}{(1+r^2)^{\frac{3}{2}}} e^{\frac{2(f-t_0)}{\k}}- 3\P ^{-2}
           \varepsilon ^2 \frac{r}{(1+r^2)^3} \Big)\\
      &=\frac{\varepsilon}{\k}e^{\frac{t_0}{\k}}<0.
 \end{aligned}
\eeQ

In hyperbolic coordinates, consider the graph $f(x)=T_0+\varepsilon e ^{-\frac{3r}{\k}}$ for certain constant $T_0 >0$ in de Sitter spacetime.
The graph is $\H$-asymptotically de Sitter ($\H=\sinh\frac{T_0}{\k}$) for sufficiently small $\varepsilon$, and
 \beQ
 \begin{aligned}
& g =\Big( \sinh^2{\frac{f}{\k}} - \frac{9\varepsilon^2}{\k^2}e ^{\frac{-6r}{\k}} \Big)dr^2
+ \k^2\sinh^2{\frac{f}{\k}}\sinh^2\frac{r}{\k}\Big(d\theta^2+\sin^2\theta \psi^2 \Big),\\
& tr _g(K)\sinh{\frac{T_0}{\k}} - \frac{3}{\k}\cosh {\frac{T_0}{\k}} =
   -\frac{12\varepsilon }{\k ^2\sinh{\frac{T_0}{\k}}}e ^{\frac{-5r}{\k}}+ O(e^{\frac{-6r}{\k}}).
 \end{aligned}
 \eeQ
Thus if $ \varepsilon <0$, then $tr _g(K)\sinh{\frac{T_0}{\k}} > \frac{3}{\k}\cosh {\frac{T_0}{\k}}$ for large $r$, and the total energy
 \beQ
 \begin{aligned}
E^{\H} _0 = &\frac{\H ^2}{16\pi}\lim_{r \rightarrow \infty}\int_{S_r}
              \mathcal{E} n^0 e^{\frac{r}{\k}} \breve{e}^2\wedge \breve{e}^3 \\
          = &\frac{\sinh^2 {\frac{T_0}{\k}}}{16\pi}\lim_{r \rightarrow \infty} \int_{S_r}
              \frac{16 \varepsilon}{\k ^2}\tanh{\frac{T_0}{2\k}}e^{-\frac{2r}{\k}} \k^2 \sinh^2 {\frac{r}{\k}}\sin{\theta}d\theta d\psi \\
          = & \varepsilon\tanh{\frac{T_0}{2\k}} \sinh^2 {\frac{T_0}{\k}}<0.
 \end{aligned}
\eeQ

The constant $\Lambda=\frac{3}{\lambda ^2}$ in the above two positive energy theorems can be chosen arbitrarily which is not necessary the cosmological
constant. Based on the experimental dates from Planck 2015\footnote{Planck 2015 results I. Overview of Scientific Results. ArXiv: 1502.01582; XIII.
Cosmological Parameters. ArXiv: 1502.01589.}, Hubble constant $H=\frac{1}{\lambda}$ is chosen and $\tilde{g} _{dS}$ is taken as
the FLRW metric with $k=0$
\beQ
\tilde g _{FLRW} = -dt ^2 + e ^{2Ht} g_\delta , \quad g_\delta=(dx ^1)^2 +(dx ^2)^2+(dx ^3)^2
\eeQ
which is isotropic and homogeneous, where
\beQ
3H^2 =\rho _m + \Lambda _c
\eeQ
and $\rho_m \cong 0.3156 \times 3H^2$ is the matter density containing dark matter, $\Lambda _c \cong 0.6844 \times 3H^2$ is the real value of
cosmological constant representing dark energy. The universe takes asymptotically FLRW metrics and satisfies the dominant energy condition,
and initial data sets are $\P$-asymptotically de Sitter. If universe's volume expansion ratio $tr_g(K)$ can be detected to be less than or equal
to $3H$ from Planck 2015, then the universe has positive total energy counting dark matter and dark energy. More precisely, the universe has more
energy than the metric $\tilde g _{FLRW}$ has. Otherwise it may have negative total energy.

The case of $\H$-asymptotically de Sitter initial data sets corresponds to FLRW metrics with $k<0$.

\mysection{$\Lambda$ <0}\ls

In this section, we assume the cosmological constant is negative. In this case spacetimes are asymptotically anti-de Sitter and initial
data sets have asymptotically hyperbolic metrics and the asymptotically zero second fundamental forms. There are also a large number of papers to devote
to define the total energy-momentum and prove its positivity in a physical manner, see, e.g. \cite{AD, HT, ACOTZ} and references therein. (It seems the
total energy was first defined in \cite{AD}, and which also contained the proof of its positivity via SUGRA, exactly as the proof for zero cosmological
constant \cite{DT}.) However, the mathematical rigorous and complete proofs were given only in \cite{Wa, CH} for asymptotically anti-de Sitter initial
data sets with the zero second fundamental form, and in \cite{M, CMT} for the initial data sets with the nonzero second fundamental form where the
energy-momentum matrix was proved to be positive semi-definite, and some energy-momentum inequalities were proved in certain specific
coordinate systems. The positive energy theorem near null infinity in \cite{Z2} gives a different energy-momentum inequality for
asymptotically anti-de Sitter initial data sets with the nonzero second fundamental form whose trace is nonpositive \cite{XZ}. In \cite{WXZ}, the
complete and rigorous study of the positive energy theorem was provided.

Denote $\kappa =\sqrt{-\frac{\Lambda}{3}}$ throughout the section. Under coordinate transformations
 \beQ
X^0=\frac{\cos(\kappa t)}{\kappa} \cosh(\kappa r),\  X^i =\frac{\sinh(\kappa r)}{\kappa} n^i,\ X^4=\frac{\sin(\kappa t)}{\kappa} \cosh(\kappa r),
 \eeQ
the anti-de Sitter metric is
\beQ
\widetilde{g}_{AdS}=-\cosh^2(\kappa r)dt^2+dr^2+\frac{\sinh^2(\kappa r)}{\kappa^2}\big(d\theta^2+\sin^2\theta d \psi^2\big)
\eeQ
and the initial data set is $(\hH ^3, \breve{g}=g _{\H}, \breve{h}=0)$. The Killing vectors $U_{\alpha\beta}$ depend on $t$ restricting on
anti-de Sitter spacetime.

Let $\breve{e} _i$, $\breve{e}^i$ and $\hconn $ be the frame, coframe and the Levi-Civita connection of the hyperbolic metric $\breve{g}$ respectively.
The initial data set $(M, g, h)$ is {\em asymptotically anti-de Sitter} if $M$ has a finite number of ends,
 \beQ
\left\{ a_{ij},\, \hconn  _k a_{ij},\, \hconn  _l \hconn _k a_{ij},\,h _{ij}, \,\hconn  _k h_{ij}\right\}=O(e^{-\kappa \tau r}), \quad \tau >\frac{3}{2},
 \eeQ
on each end, where $a_{ij}=g(\breve{e} _i,\breve{e} _j)-\breve{g}(\breve{e}_i,\breve{e}_j), h_{ij}= h (\breve{e} _i,\breve{e} _j )$. Moreover, let $R$, $\nabla$,
$\rho_z$ be the scalar curvature and the Levi-Civita connection of $g$, and the distance function with respect to $z \in M$ respectively,
 \beQ
(R +6 \kappa ^2 )e^{\kappa \rho_z},\quad (\nabla^j h_{ij}-\nabla_i tr_{g}(h))e^{\kappa \rho_z} \in L^1(M).
 \eeQ

In 1985, Henneaux and Teitelboim defined the total energy-momentum $J_{ab} ^{HT}$ for asymptotically AdS spacetimes associated to $U_{ab}$ \cite{HT}.
In our notation, these quantities are
 \beQ
\begin{aligned}
E_0=&\frac{\kappa}{16\pi}\int_{S_\infty}\mathcal{E} U_{40}^{(0)}\breve{\omega},\\
c_{i}(t)=&\frac{\kappa}{16\pi}\int_{S_\infty}\mathcal{E} U_{i4}^{(0)}\breve{\omega}
       +\frac{\kappa}{8\pi}\int_{S_\infty}\mathcal{P}_{A}U_{i4}^{(A)} \breve{\omega},\\
c'_{i}(t)=&\frac{\kappa}{16\pi}\int_{S_\infty}\mathcal{E} U_{i0}^{(0)}\breve{\omega}
        +\frac{\kappa}{8\pi}\int_{S_\infty}\mathcal{P}_{A}U_{i0}^{(A)} \breve{\omega},\\
J_{i}=&\frac{\kappa}{8\pi}\int_{S_\infty}\mathcal{P}_{A}V_{i}^{(A)} \breve{\omega},
\end{aligned}
 \eeQ
where $\mathcal{E} =\hconn ^{j}g_{1j}- \hconn  _1 tr _{\breve{ g}}(g ) +\kappa (a_{22}+a_{33})$, $\mathcal{P}_{j} =h_{j1}-g_{j1}tr_{\breve{g}}(h)$,
$\breve{\omega}= \breve{e}^2\wedge \breve{e}^3$, $U_{\alpha\beta}=U_{\alpha\beta}^{(\gamma)}\breve{e}_{\gamma}$,
$V_i=\frac{1}{2}\varepsilon _{ijk} U_{jk}=V_{i}^{(A)} \breve{e}_{A}$, $A=2,3.$

Denote ${\bf c} =(c_1, c_2, c_3)$, ${\bf c}' =(c' _1, c' _2, c' _3)$, ${\bf J} =(J_1, J_2, J_3)$ and
\beQ
\begin{aligned}
L=&\left(|{\bf c}|^2 +|{\bf c}'|^2 +|{\bf J}|^2 \right)^{\frac{1}{2}},\\
A=&\left(|{\bf c} \times {\bf c}' |^2+|{\bf c} \times {\bf J}|^2 +| {\bf c}'\times {\bf J} |^2\right)^{\frac{1}{4}},\\
V=&\left(\varepsilon_{ijl}c_i c_{j}' J_l \right) ^{\frac{1}{3}}.
\end{aligned}
\eeQ
The quantities $2L$, $2A^2$ and $V^3$ are the (normalized) length, surface area and volume of the parallelepiped spanned by ${\bf c}$, ${\bf c}'$ and
${\bf J}$, which are independent on $t$. Clearly, $L^2 \geq 3V^2$, and $|{\bf c}|^2 + |{\bf c}' |^2$ is independent on $t$.

The pseudo-Euclidean space $\R^{3,2}$ has two timelike Killing vectors and three spacelike Killing vectors.
Physically, $E$ measures the rotation on the plane $(X^0, X^4)$, $c_i$ measures the rotation on the plane $(X^i, X^4)$,
$c' _i$ measures the rotation on the plane $(X^0, X^i)$ and $J_i$ measures the rotation on the plane $(X^j, X^k)$ where $\{i,j,k\}$ is the even
permutation of $\{1,2,3\}$. But these rotations are all observed from a curved hyperboloid, so they contain both
translation and rotation of an asymptotically anti-de Sitter spacetime. This indicates that we can not simply refer them as the center of mass
as well as the total angular momentum. The total effect of translation and rotation is given by the parallelepiped spanned by
${\bf c}$, ${\bf c}'$ and ${\bf J}$ which can be measured from its length of the edges, surface area and the volume.

In 2001, Wang established the mathematically rigorous and complete proofs for asymptotically anti-de Sitter initial data sets with $a_{12}=a_{13}=0$
and the zero second fundamental form \cite{Wa}. And these extra conditions on $a$ were removed by Chr\'{u}sciel-Herzlich later \cite{CH}. Under the
dominant energy condition, they proved, at $t=0$,\\
{\em (i) $E_0 \geq \sqrt{c_1 ^2 + c_2 ^2 +c_3 ^2}$ for each end;\\
(ii) That $E=0$ for some end implies $M$ is hyperbolic.}

In 2006, Maerten used the Lorentzian setting to treat $M$ as a spacelike hypersurface in spacetime ${\bf L} ^{3,1}$. Using Witten's approach and
studying the restriction of spin geometry of ${\bf L} ^{3,1}$ over $M$ with the nonzero second fundamental form, he proved that there is a positive
semi-definite $4\times 4$ Hermitian matrix ${\bf Q}$ involving $E_0$, $c_i(0)$, $c_i '(0)$ and $J_i$ under the dominant energy condition \cite{M}.
Soon later, Chr\'{u}sciel, Maerten and Tod \cite{CMT} proved, at t=0, if $E_0 > \sqrt{c_{1} ^2 + c_{2} ^2 +c_{3} ^2} $, one can make $SO(3,1)$
coordinate transformations such that
 \beQ
\sqrt{E_0 ^2 - c_{1} ^2 - c_{2} ^2 -c_{3} ^2} \rightarrow \bar E_0,\ c_i, c'_2, J_1, J_2 \rightarrow 0,\
c' _1 \rightarrow \bar c _1 ',\ c' _2 \rightarrow \bar c _2 ', \ J_3 \rightarrow \bar {J _3}.
 \eeQ
In the new coordinates, which we refer to the ``center of AdS mass" coordinates, they proved\\
{\em (i) $\bar E_{0}\geq \sqrt{| \bf{\bar c '}|^2+| \bf{\bar J} |^2+ 2| \bf{\bar c '}\times \bf{\bar J} |}$ for each end;\\
(ii) That $E=0$ for some end implies $M$ has only one end, and ${\bf L} ^{3,1}$ is anti-de Sitter along $M$.}

In \cite{WXZ}, it was shown that the trace, sum of the second-order minors, sum of the third-order minors and the determinant of ${\bf Q}$ are
 \beQ
 \begin{aligned}
tr{\bf Q}=&4E_0, \quad {\bf Q} ^{(2)}=6E_0 ^2 -2L^2, \quad
{\bf Q} ^{(3)}=4E_0 (E_0^2-L^2)+8V^3, \\
\det {\bf Q} =&(E_0^2-L^2 )^2+8 E_0 V^3-4A^4.
 \end{aligned}
 \eeQ
But Witten's argument indicates only that these quantities are nonnegative, which do not yield $E_0 > \sqrt{c_{1} ^2 + c_{2} ^2 +c_{3} ^2 } $ in general
when the second fundamental form is nonzero, required for the center of AdS mass coordinate transformations. (e.g.
${\bf Q} ^{(2)} > 0 \Longrightarrow E_0 >\sqrt{\frac{1}{3}(|{\bf c}|^2 +|{\bf c}'|^2 +|{\bf J}|^2)}
\not \Longrightarrow E_0 > \sqrt{c_{1} ^2 + c_{2} ^2 +c_{3} ^2}$.)
On the other hand, the form of Chr\'{u}sciel-Maerten-Tod's energy-momentum inequality is not $SO(3,1)$ invariant. It changes when it is transformed
back to the non-center of AdS mass coordinates. These problems are the motivation to establish the inequality for Henneaux and Teitelboim's total energy-momentum
in general non-center of AdS mass coordinates at general $t$, and the positive energy theorem in this case was proved by Wang, Xie and Zhang \cite{WXZ}.

{\em The positive energy theorem: If $(M, g, h)$ is an asymptotically anti-de Sitter initial data set in spacetime ${\bf{L}} ^{3,1}$ with negative cosmological constant $\Lambda <0$,
with possibly a finite number of apparent horizons which are inner boundaries of $M$, each of them is topological $S^2$ whose mean curvature $H$
satisfies $H \pm tr_g (h| _{S^2}) =0$. Suppose the dominant energy condition holds, then\\
(i) $E_0 \geq \sqrt{L^2 -2V^2 +2 \big(\max\{A^4 -L^2 V^2, 0\}\big)^\frac{1}{2} }$ for each end;\\
(ii) That $E=0$ for some end implies $M$ has only one end, and ${\bf L} ^{3,1}$ is anti-de Sitter along $M$.}

If three vectors ${\bf c}$, ${\bf c}'$, ${\bf J}$ are linearly dependent, i.e, $V=0$, then the above energy-momentum inequality reduces to
$E_0 \geq \sqrt{L^2 +2 A^2}$ which provides more general energy-momentum inequality than Chr\'{u}sciel-Maerten-Tod's.

In \cite{WXZ}, Wang, Xie and Zhang also proved that $\det{\bf Q}$ is invariant under the following admissible coordinate transformation on ends
 \beQ
\begin{aligned}
\hat t = t +o(e^{-\frac{5\kappa}{2}r}), \quad & \breve{e} _0 (\hat t)=\breve{e} _0 {(t)} +o(e^{-\frac{7\kappa}{2}r}),\\
\hat r = r +o(e^{-\frac{3\kappa}{2}r}), \quad & \breve{e} _1 (\hat r)=\breve{e} _1 {(r)} +o(e^{-\frac{3\kappa}{2}r}),\\
\hat \theta ^A = \theta ^A +o(e^{-\frac{5\kappa}{2}r}), \quad & \breve{e} _B (\hat \theta ^A)=\breve{e} _B {(\theta ^A)} +o(e^{-\frac{7\kappa}{2}r}).
\end{aligned}
 \eeQ
Thus it serves as the geometric invariant of asymptotically anti-de Sitter spacetimes. Moreover,
$\det {\bf Q}=\Big(\frac{\kappa}{16\pi}\Big)^4 \big(I_1^2+2I_2 \big)$, where
\beQ
I_1 =\frac{1}{2}J^{HT} _{ab}J^{HT, ab},\quad
I_2 =\frac{1}{2}J _{a}^{HT, b}J _{b}^{HT, c}J_{c}^{HT, d}J_{d}^{HT, a}-\frac{1}{4}(J_{a}^{HT, b}J_{b}^{HT,a})^2.
\eeQ
are two $O(3,2)$ Casimir invariants \cite{HT}. Therefore $\sqrt[4]{\det{\bf Q}}$ is $O(3,2)$ invariance and serves as the total rest mass of
asymptotically anti-de Sitter spacetimes.

Finally, we would like to remark that the total energy-momentum for negative cosmological constant may have some physical implications in the theory
of strongly coupled superconductors based on the point of view of Anti-de Sitter/Conformal Field Theory correspondence. We refer to the recent
introductory overview of relevant theory \cite{CLLY} and references therein.

\bigskip

{\footnotesize {\it Acknowledgement.}
This work is partially supported by the National Science Foundation of China (grant 11171328) and the project of mathematics and interdisciplinary
sciences of Chinese Academy of Sciences.}

\end{document}